\documentclass[12pt]{article}

\usepackage{latexsym}
\usepackage{amssymb}
\usepackage{theorem}
\usepackage{amsmath}
\usepackage{amscd}

\usepackage[pdftex]{graphicx}
\usepackage{ulem}

\newtheorem{thm}{Theorem}[section]
\newtheorem{cor}[thm]{Corollary} 
\newtheorem{lem}[thm]{Lemma} 
 
\newtheorem{rem}[thm]{Remark} 
\newtheorem{exam}[thm]{Example} 
 
\newtheorem{defin}[thm]{Definition} 
 
\newtheorem{fact}[thm]{Fact}

\newcommand{\enD}{\hfill $\Box$\vspace{3truemm} \par}
\newcommand{\R}{\mathbb{R}}

\def\proof{\noindent {\sl Proof} :\;  }

\begin{document}
\title{Visualization of curvature on curve and surface by tangential angle parametrization}

\author{Yutaro Kabata, Shigeki Matsutani and Yuta Ogata}

\date{\today}

\maketitle
\begin{abstract}
We propose a unified method to visualize curvature on planar curves and surfaces of revolution using the tangential angle parameter. For plane curves, placing markers at equal increments of the tangential angle reveals local bending features and naturally highlights inflection points and vertices. This approach extends to surfaces of revolution, where curvature lines drawn at equal tangential angle steps reflect principal curvature variations and naturally expose ridge and parabolic curves. Our method provides clear, consistent visualizations without arbitrary parameter tuning, offering geometric insight for both analysis and design applications.
\end{abstract}

\renewcommand{\thefootnote}{\fnsymbol{footnote}}





\section{Introduction}
Curvature is a fundamental geometric concept that quantifies how curves and surfaces bend. Clear and consistent visualization of curvature is essential in both theoretical research and practical applications. In this paper, we demonstrate a simple yet powerful method for visualizing curvature on planar curves and surfaces of revolution by using the curve's {\it tangential angle}~$\theta$.


Visualizing curvature can be hard when a shape contains local parts with complicated
details in large regions.  For example, engineers often need to assess curvature of
a small feature without losing sight of its meaning on the overall context.  Meeting
this demand requires a visual representation that preserves both local detail and
global structure.  In practice, this means projecting the bending behaviour of curves
and surfaces onto a two-dimensional plane in a perceptually faithful way.  The
computer-graphics (CG) and computer-vision (CV) communities have long addressed this issue \cite{Cipolla-Giblin, Damon-Giblin-Haslinger,Hertzmann2000, Dong2006, Palacios2007,  Koenderink, Koenderink1992, KD1998}.

On planar curves, sampling at equal arc-length steps may miss sharp bends or waste
markers on nearly straight stretches.  On surfaces, particularly in computer graphics
and manufacturing, curvature lines are commonly used to illustrate surface bending,
but their spacing is usually chosen arbitrarily \cite{Iarussi2015,Takezawa2019}.
Fixed intervals can lead to uneven coverage—some areas become overcrowded with
lines while others are sparsely detailed—because line density depends on both
curvature magnitude and the chosen parameter.  A parameterization that
automatically adjusts line placement to local curvature would remove this
arbitrariness, enabling consistent, informative visualizations without manual tuning.

\begin{figure}
\centering{
\includegraphics[width=120mm]{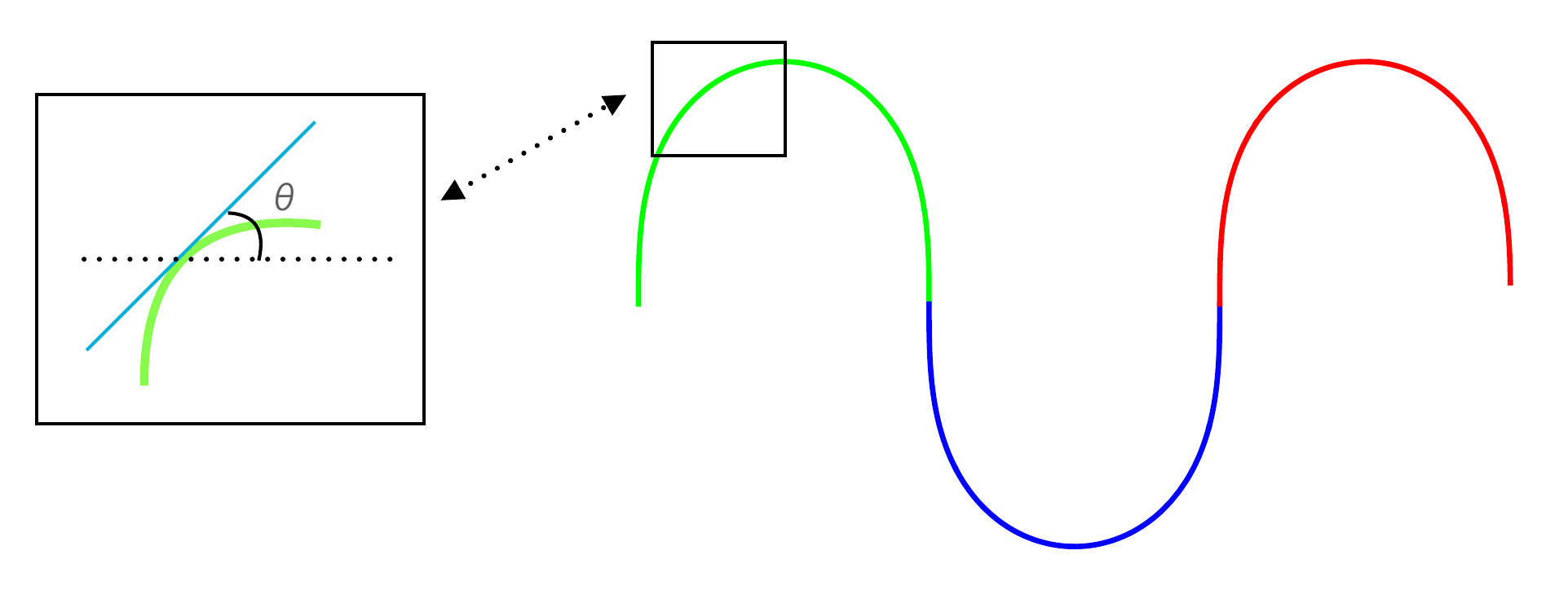}
}
\caption{Tangential angle \(\theta\) as the angle formed by the tangent vector of the curve and the positive~\(x\)-axis. Inflection points split the curve into colored segments, and the tangential angle \(\theta\) is taken separately on each segment.}
\label{fig:tangentialangle}
\end{figure}

\begin{figure}
{\centering
\includegraphics[width=80mm]{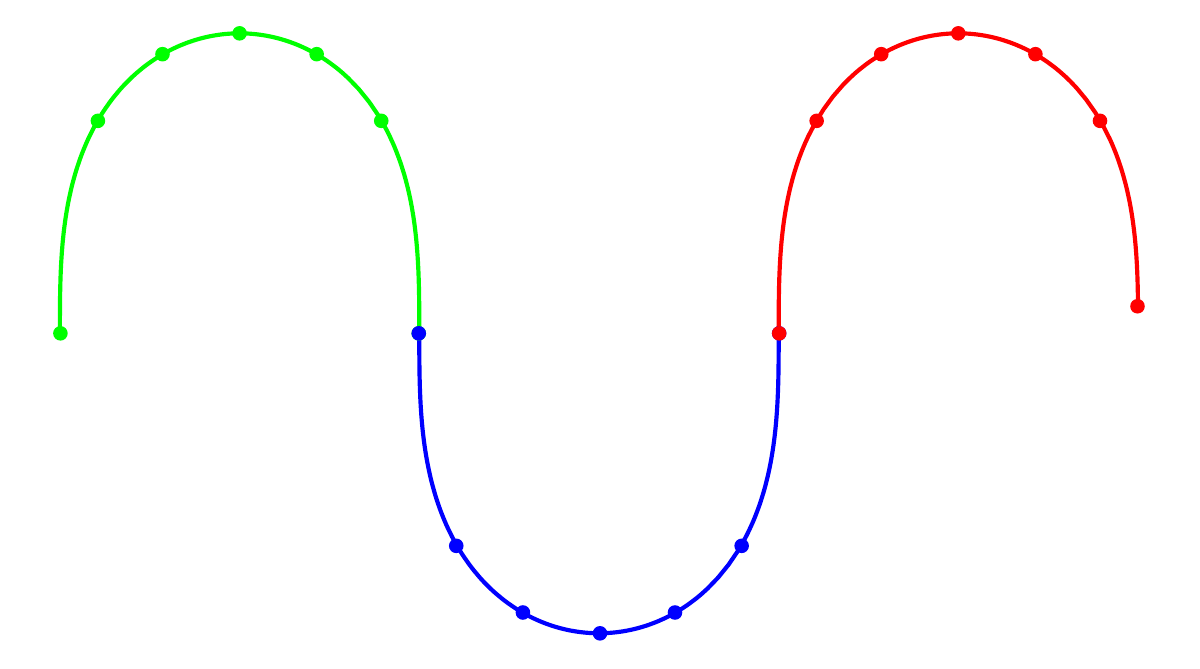}
\hfill
\includegraphics[width=80mm]{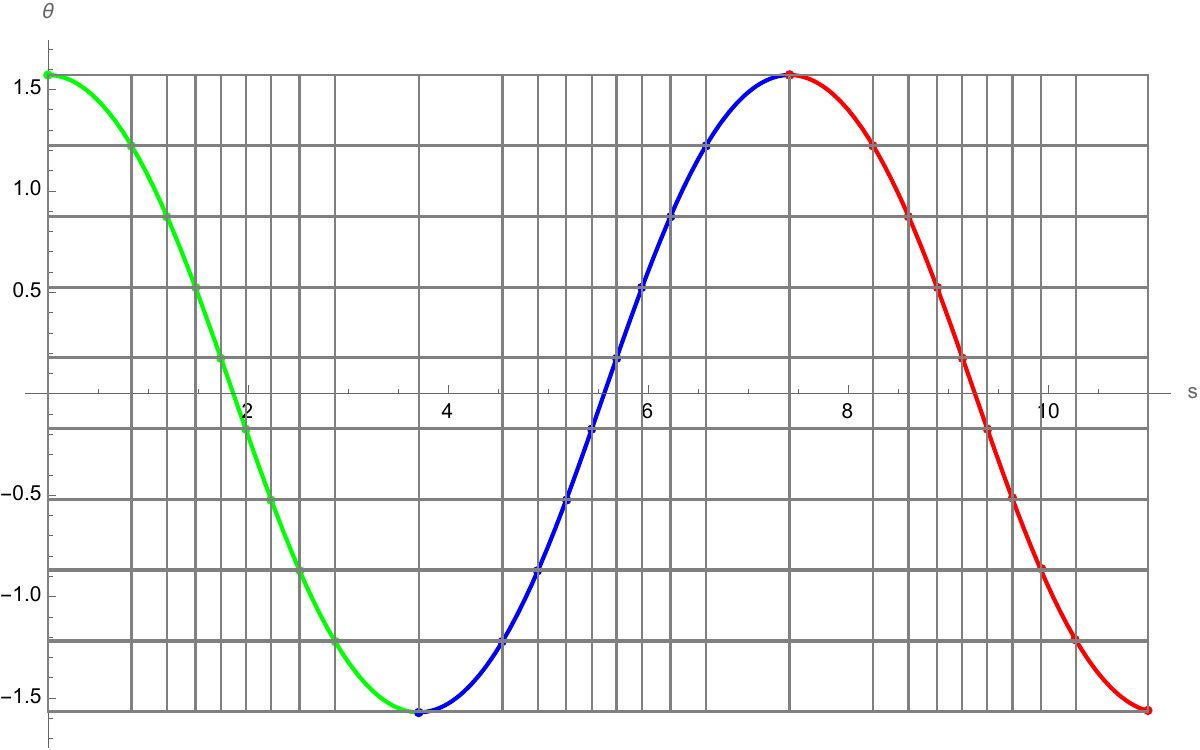}
}
\caption{Left: An elastica curve split at inflection points (\(\kappa=0\)) with markers placed at equal increments of the tangential angle \(\theta(s)=\int_{c}^{s}\kappa(u)\,\mathrm{d}u\). Right: Plot of \(\theta(s)\) versus arc length \(s\), showing that equal steps in \(\theta\) correspond to non-uniform steps in \(s\).}
\label{fig:longelastica}
\end{figure}

Our approach begins with planar curves. We define the {\it tangential angle} $\theta(s)=\int_{c}^{s}\kappa(u)\,\mathrm{d}u$, where $\kappa(s)$ is the curvature, $s$ the arc-length parameter and $c$ a real constant. See Figure \ref{fig:tangentialangle}. We place markers along the curve at points where $\theta$ increases by a fixed step. This method makes it easy to see how curvature changes along a curve. In the left side of Figure~\ref{fig:longelastica}, we cut the elastica at its inflection points (where $\kappa=0$) and place markers so that the tangential angle increases by the same amount each time. The markers get closer together where the curve bends more. The right side of Figure~\ref{fig:longelastica} plots $\theta(s)$ against the arc length $s$: equal steps in $\theta$ correspond to uneven steps in $s$, highlighting that the tangential angle functions as an independent parameter, separate from arc length. These observations are mathematically verified by Theorem \ref{thm:curvthm}.

\begin{figure}[ht]
  \centering
  \includegraphics[width=0.7\linewidth]{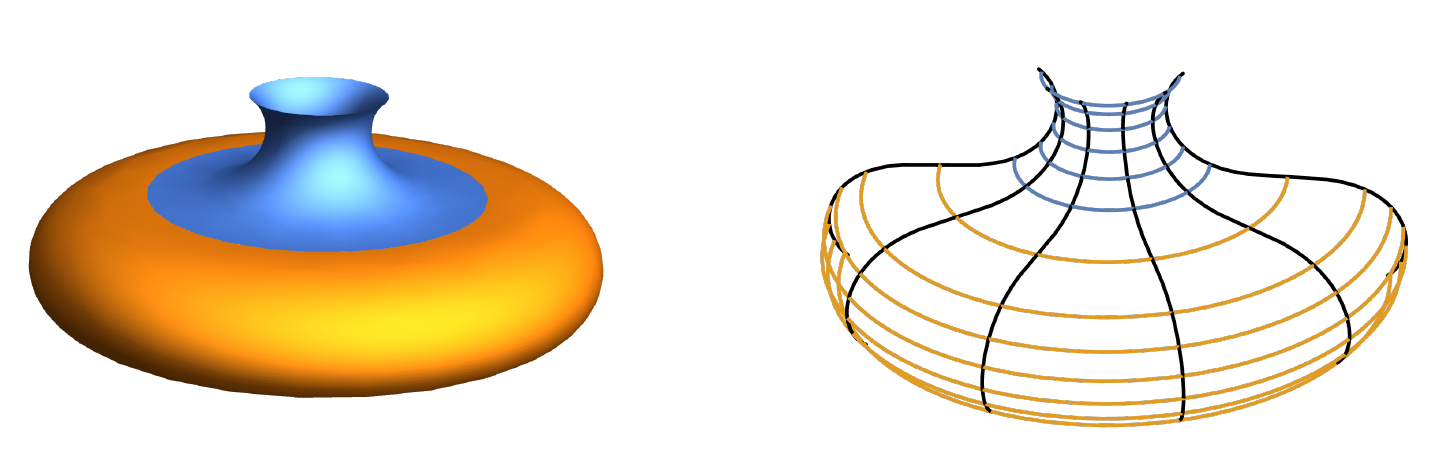}
  \caption{Surface of revolution with the profiles being an Euler spiral. The left figure shows the half part drawn by only curvature lines: especially, parallel circles drawn at equal $\theta$ intervals along the profile curve. Here, the parabolic curve divides the surface into two regions: the elliptic region colored with orange and the hyperbolic region colored with blue.}
  \label{fig:revolutionintro}
\end{figure}

Next, we extend this idea to surfaces of revolution. As in Figure \ref{fig:revolutionintro}, along the profile curve, we draw {parallel circles} at equal $\theta$ intervals, where $\theta$ is a tangential angle parameter of the profile curve (plane curve). These circles are curvature lines on the surface: dense bands indicate regions of high principal curvature, while sparse bands indicate low curvature (Figure \ref{fig:revolutionintro}). This observation fact is mathematically verified by Corollaly \ref{cor:surfaceofrevo}. Furthermore, our method gives a natural segmentation of a surface by parabolic curves and ridge curve. Unlike conventional visualizations that rely on shading, lighting effects, or reflections, our approach conveys geometric structure purely through curvature lines.

This unified framework enhances our geometric insight to understand curves and surfaces from their graphics. It also helps guide intuitive design and manufacturing decisions, enables efficient surface segmentation, and removes the arbitrary parameter choices of traditional curvature-line displays—ensuring clear, consistent visualizations.

In Section 2, we define the tangential angle parameter for planar curves and demonstrate its effectiveness in visualizing curvature. We show that using this parameter naturally stratifies the curve by vertices and inflection points, and we provide several examples illustrating how curvature affects the distribution of plotted points. In Section 3, we extend this idea to surfaces of revolution. We demonstrate that the tangential angle parameter induces a meaningful coordinate system that highlights geometric features such as ridge and parabolic curves. This framework not only unifies the visualization of curvature across dimensions but also enables consistent and interpretable curvature-based graphics without relying on arbitrary parameter choices.

\bigskip
\noindent
{\bf Acknowledgement}. 
This work is partially supported by JSPS KAKENHI Grant Number JP 20K14312, 24K22308 and 24K06722.

\section{Tangential angle of plane curve}
Let $\R^2$ be the Euclidian plane equipped with the standard inner product $<\cdot,\cdot>$ and Euclidian norm $|\cdot|$.
Let $\gamma\colon I\to \R^2$ be a smooth ($C^\infty$) non-singular curve with an open interval $I$.
The parameter $s$ of $\gamma$ is called {\it the arc-length parameter of $\gamma$}, if $|\frac{d\gamma}{ds}|=1$ holds.
Putting $e=\frac{d\gamma}{ds}$, then {\it the curvature $\kappa$ of $\gamma$} is defined as the function such that
$$
\frac{de}{ds}=\kappa \nu, 
$$
where $\nu\colon I\to S$ is the Gauss map.
It is well-known that $\frac{1}{\kappa}$ is the radius of the circle which approximates the curve locally. 
\begin{defin}\label{def:curveintpara}
Take $J\subset I$ as an open interval on which $\kappa\not=0$ holds, and define 
$$
\theta\colon J \to \theta(J),\;\;\theta(s)=\int^s_c \kappa(s) d s,
$$
with a constant $c\in J$. The parameter $\theta$ is called {the tangential angle of $\gamma$}.
\end{defin}
Note that $\theta$ is a diffeomorphism since $\theta'=\kappa\not=0$ on $J$.

\begin{fact}[cf. Proposition 1.20 in \cite{Porteous}]
$\theta$ is regarded as an angle between the tangent vector $e$ and some axis in the plane.
\end{fact}

\begin{rem}
{\rm In \cite{Porteous}, the parameter $\theta$ is called {\it the angle of contingence}; and
the plane curve $\gamma$ parametrized by $\theta$ is called {\it a unit-angular velocity curve}.
}
\end{rem}


Note that the point $\gamma(t_0)$ with $\kappa(t_0)=0$ is called {\it an inflection}, where the tangent line have a contact of the order of more than $2$. Thus we can take the tangential angle as a parameter of $\gamma$, after we divide $\gamma$ by eliminating inflections.

\begin{lem}[cf. Theorem 1.21 in \cite{Porteous}]\label{lem:thetaconst}
Given a nonzero function $\kappa\colon J\to\R$. Then the curvature of the following curve is $\kappa$:
\begin{equation}\label{eq:construction}
\gamma(\theta)=\int^{\theta}_c\frac{1}{\kappa({\theta})}(\cos {\theta},\sin{\theta})d{\theta}.
\end{equation}
Here $c\in J$ is a constant.
\end{lem}
\proof
Recall that the plane curve $\gamma$ with respect to the arc-length parameter $s$ is constructed from a given curvature $\kappa$ as follows: 
\begin{equation*}
\gamma(s)=\int^{s}_{\tilde{c}}(\cos \theta(s),\sin\theta(s))ds
\end{equation*}
for a suitable real constant $\tilde{c}$.
\enD

The following is a main result of the present paper showing that the tangential angle $\theta$ is useful to efficiently visualize information of the curvature.
\begin{thm}\label{thm:curvthm}
Suppose that $\gamma\colon J \to \R^2$ doesn't have inflections.
Then, we have
\begin{equation}\label{eq:mainthm}
\frac{d }{d\theta}\left|\frac{d \gamma}{d\theta}\right|^2=
-\frac{2}{\kappa^4}\frac{d\kappa}{ds}.
\end{equation}
In particular, $\frac{d }{d\theta}\left|\frac{d \gamma}{d\theta}\right|^2(t_0)=0$ if and only if $\gamma(t_0)$ is the vertex of $\gamma$.
\label{curveparametrizationthm}
\end{thm}
\proof
Taking a suitable motion on the plane, we can assume that a given curve $\gamma$ is expressed as the form of (\ref{eq:construction}).
Hence, we have
$$
\left|\frac{d \gamma}{d\theta}\right|^2=\frac{1}{\kappa^2},
$$
and
$$
\frac{d \kappa}{d\theta}=\frac{d \kappa}{ds}\cdot\frac{d s}{d\theta}=\frac{d \kappa}{ds}\cdot\frac{1}{\kappa} .
$$
Thus we have the statement by direct calculations:
$$
\frac{d }{d\theta}\left|\frac{d \gamma}{d\theta}\right|^2=\frac{d }{d\theta}\left(\frac{1}{\kappa^2} \right)=-\frac{2}{\kappa^3}\cdot\frac{d \kappa}{d\theta}=
-\frac{2}{\kappa^4}\frac{d\kappa}{ds}.
$$
\enD

\begin{rem}{\rm
Notice that Definition \ref{def:curveintpara} and Theorem  \ref{thm:curvthm} are easily extended to the case where the curve is embedded in a general Euclidian $n$-space.}
\end{rem}

From Theorem \ref{curveparametrizationthm}, we observe the following: For a curve parameterized by the new variable $\theta$, the distribution of points spaced evenly along the parameter intervals $I$ corresponds to the distribution of curvature. The absolute value of the right-hand side of (\ref{eq:mainthm}) decreases as $|\kappa|$ increases. Notably, (\ref{eq:mainthm}) is proportional to the derivative of the curvature $\frac{d\kappa}{ds}$, that is, the sign of $\frac{d }{d\theta}\left|\frac{d \gamma}{d\theta}\right|^2$ is determined by the sign of $\frac{d\kappa}{ds}$.


We show the above property with several examples of plane curves.

\begin{figure}
\centering{
\includegraphics[width=100mm]{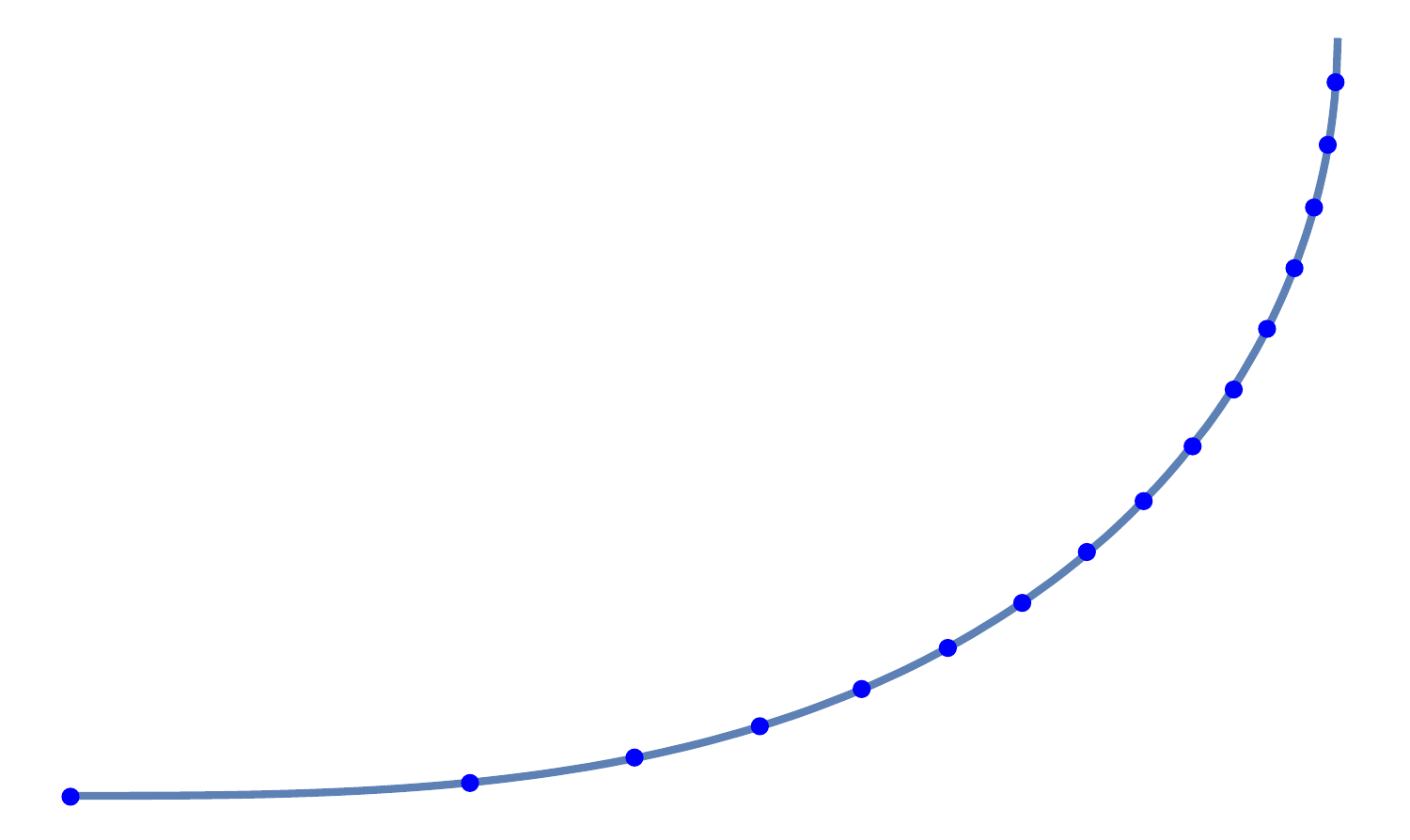}
}
\caption{An elastica. The dots are plotted according to the parameter $\theta$.}
\label{fig:elastica}
\end{figure}

\begin{exam}
Consider an elastica $\gamma\colon (0,1) \to \R^2,\; x\mapsto(x,\int^x_0\frac{x^2dx}{\sqrt{1-x^4}})$. It is well-known that the arc-length is expressed as the lemniscate integral $s(x)=\int^x_0\frac{dx}{\sqrt{1-x^4}}$ and the inverse function $x(s)$ is the lemniscate sine.
 Here, the curvature is expressed as $\kappa=2x$, and thus $\theta=\int^s_02x(s)ds=\int^{x}_0\frac{2xdx}{\sqrt{1-x^4}}=\arcsin (x^2)$. Thus we have $x=\pm \sqrt{\sin \theta}$.
See Figure \ref{fig:elastica}. Here the markers in each curve are plotted with respect to the parameter $\theta$, which are plotted in equidistance in the domains $\R_{>0}$ and $\R_{<0}$. It is clearly seen that the intervals of the plotted markers get narrower (resp. wider) as $|\kappa|$ becomes bigger (resp. smaller). 
\end{exam}

\begin{figure}
\centering{
\includegraphics[width=110mm]{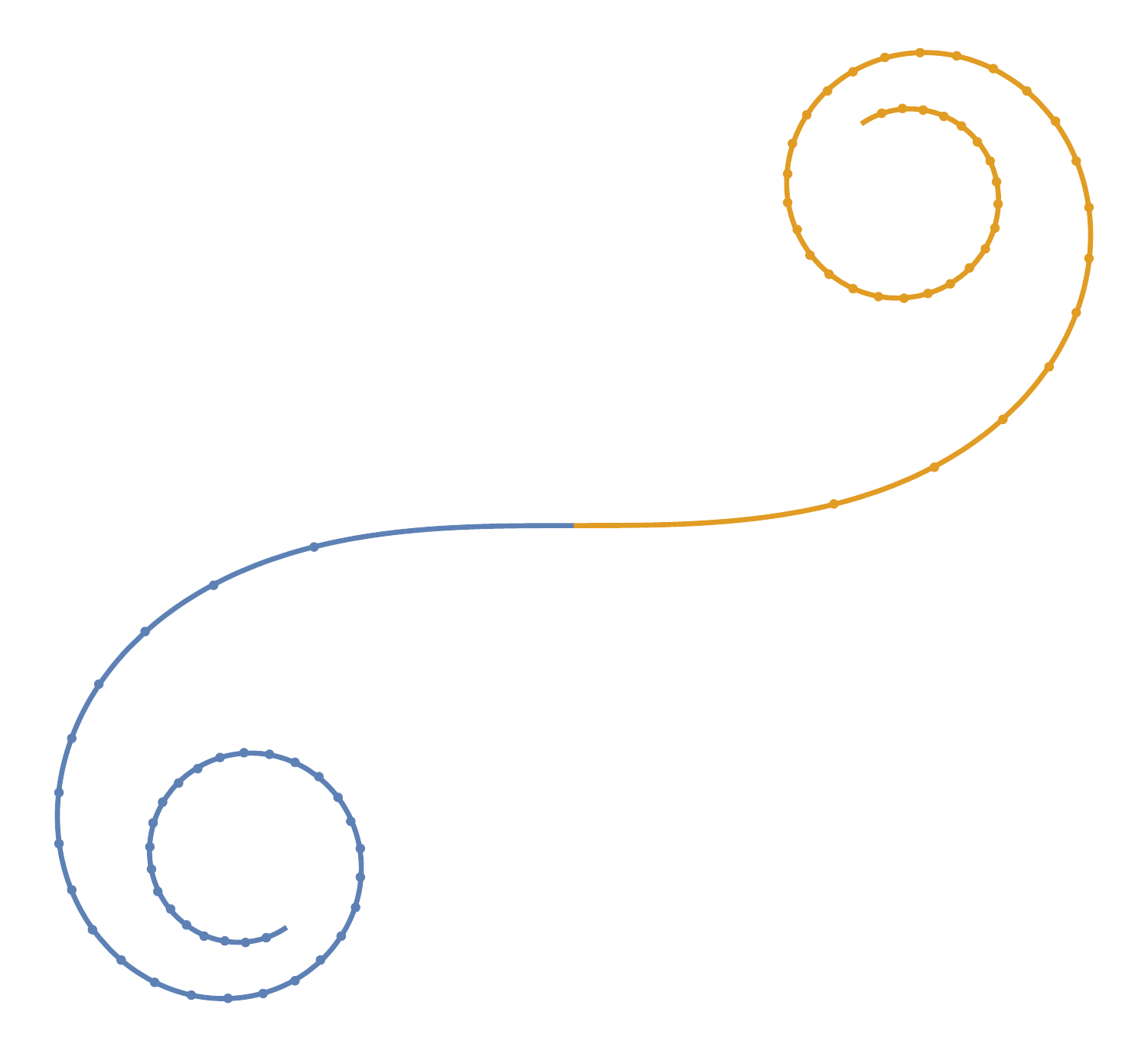}
}
\caption{The Euler spiral with $\kappa(s)=s$. The dots are plotted according to the parameter $\theta$. The curve is divided by the inflection (the point with $\kappa=0$).}
\label{fig:kappa=s}
\end{figure}

\begin{exam}
Consider a plane curve $\gamma\colon \R \to \R^2,\; s\mapsto(\int_0^s\cos\frac{t^2}{2}dt,\int_0^s\sin\frac{t^2}{2}dt)$ with the arc-length parameter $s$, which is called an Euler spiral (also called a clothoid or a Cornu spiral). Here, the curvature is expressed as $\kappa(s)=s$. Since $\kappa$ takes $0$ at $s=0$, we need divide the domain $\R$ into $\R_{>0}$ and $\R_{<0}$ in order to take the parameter $\theta$ for this case. In this setting, $\theta$ is given as $\theta=\frac12 s^2$, i.e. $s=\pm \sqrt2 \sqrt{\theta}$ for $\R_{>0}$ and $\R_{<0}$. See Figure \ref{fig:kappa=s}. Here the domain $\R_{>0}$ corresponds to the part with the color orange and $\R_{<0}$ corresponds to the part with the color blue. In addition, the markers in each curve are plotted with respect to the parameter $\theta$, which are plotted in equidistance in the domains $\R_{>0}$ and $\R_{<0}$. It is clearly seen that the intervals of the plotted markers get narrower (resp. wider) as $|\kappa|$ becomes bigger (resp. smaller). 
\end{exam}

\begin{figure}
\centering{
\includegraphics[width=110mm]{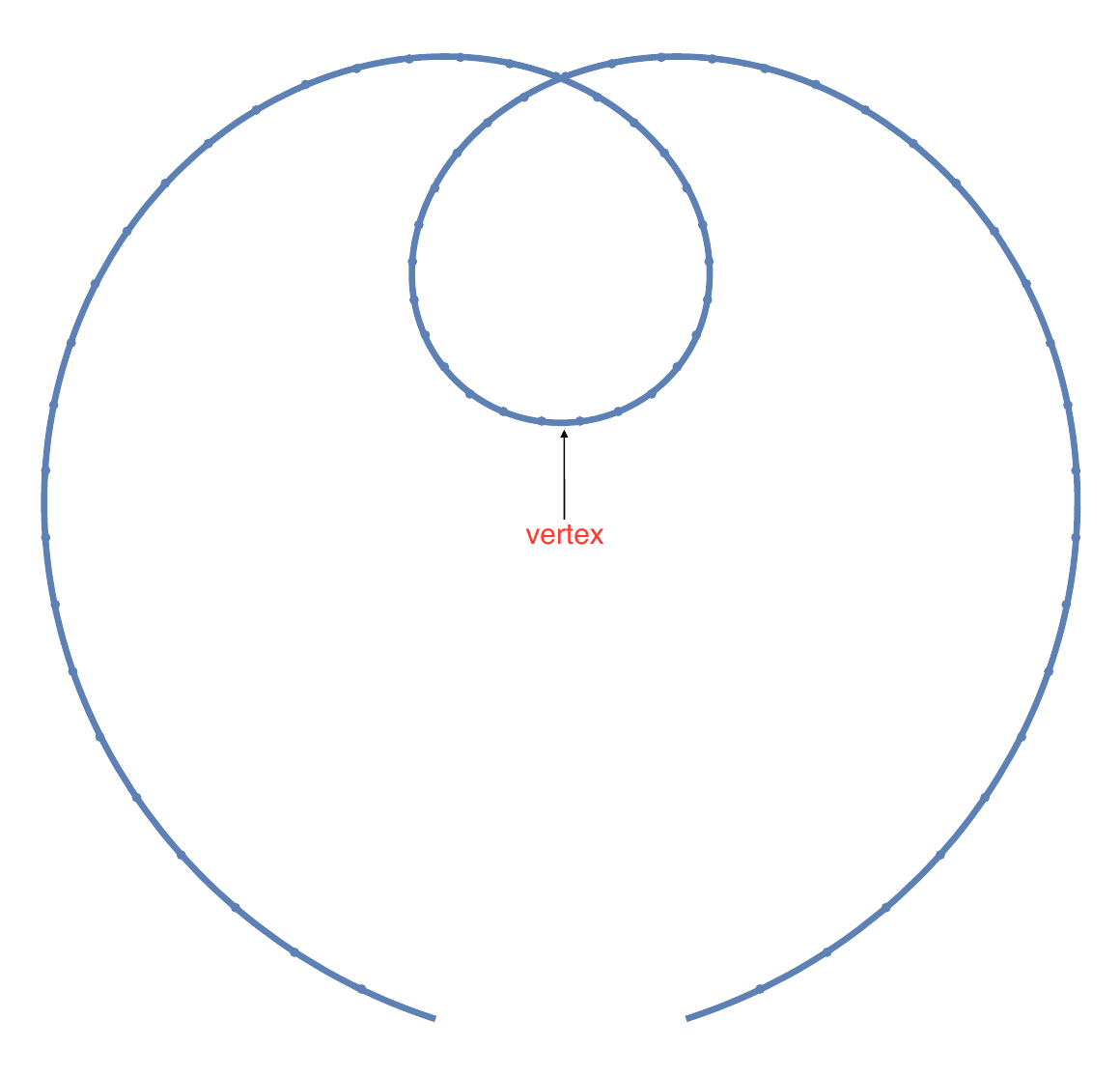}
}
\caption{The curve with $\kappa(s)=1+s^2$. The markers are plotted according to the parameter $\theta$.}
\label{fig:k=s2+1}
\end{figure}

\begin{exam}\label{ex:curvewithvertex}
Consider a plane curve $\gamma\colon \R\to \R^2$ constructed from the curvature $\kappa(s)=1+s^2$ with an arc-length parameter $s$. Setting $\theta(s)=\int^s_0 \kappa(s) ds=s+\frac13s^3$, from Lemma \ref{lem:thetaconst}
 we can construct a curve 
$$
\gamma(\theta)=-\frac{\sqrt[3]{2} \left(\sqrt{9 \theta ^2+4}-3 \theta \right)^{2/3}-2 }{2^{2/3} \sqrt[3]{\sqrt{9 \theta ^2+4}-3 \theta }}
(\cos \theta, \sin \theta).
$$
Since $\frac{d\kappa}{ds}(0)=0$, $c(0)$ is a vertex. See Figure \ref{fig:k=s2+1}. We notice that the interval of the markers are the narrowest at the vertex.
\end{exam}

\begin{rem}{\rm
In the above examples, we see the following important properties of our parametization:
\begin{enumerate}
\item With the parameter $\theta$, the information of how the curve is locally bended is efficiently visualized by the markers plotted along $\theta$: the intervals of markers is related to the differential of the curvature (see Theorem \ref{curveparametrizationthm}).
\item
In the process of taking the parameter $\theta$, we naturally get a stratification of the curve by vertices ($\kappa'=0$) and inflection ($\kappa=0$). This suggests that a suitable stratification of geometrical objects gives efficient understanding of it.
\end{enumerate} 
}
\end{rem}

\begin{figure}
\centering{
\includegraphics[width=100mm]{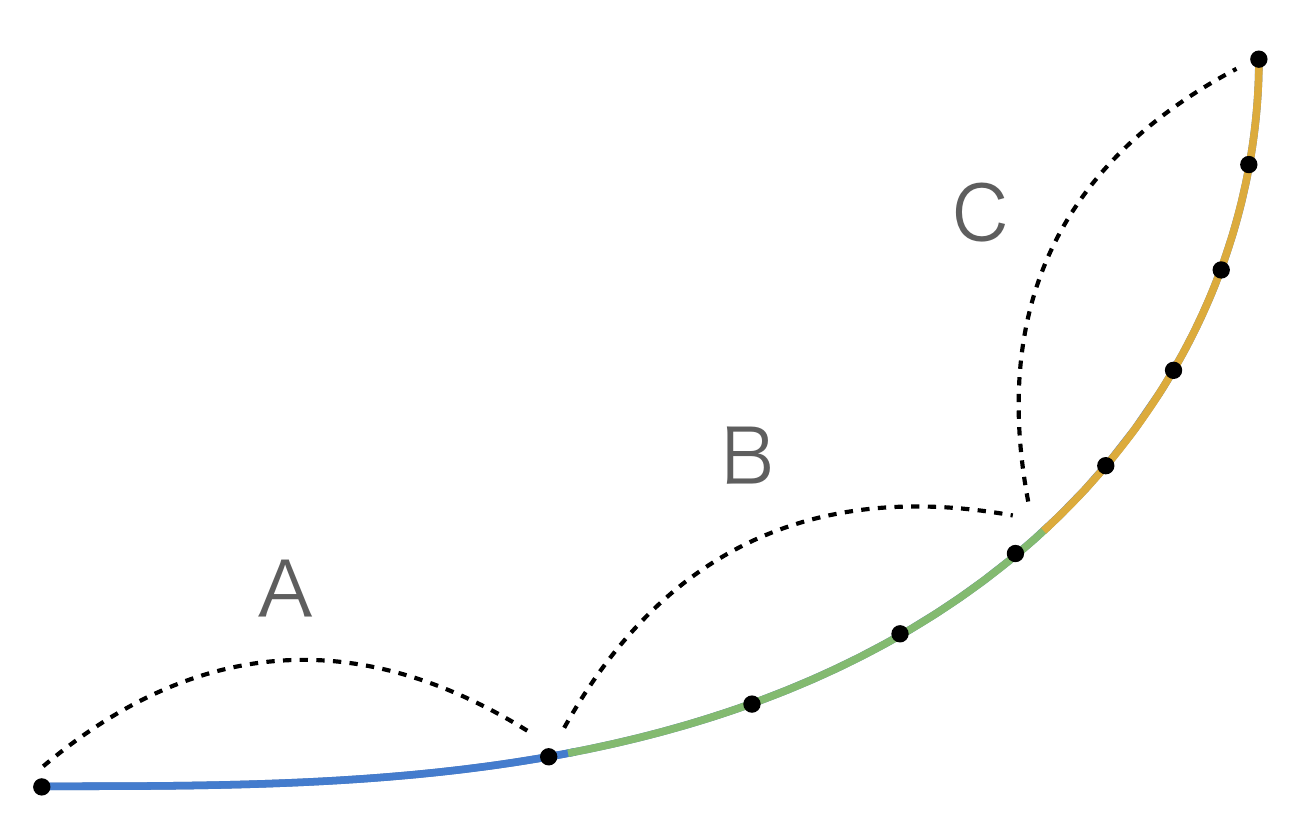}
}
\caption{The curve is divided into three equal parts in terms of arc length for each segment \( A \), \( B \), and \( C \), and the dots are evenly spaced with respect to \( \theta \).
}
\label{fig:elasticasegmentation}
\end{figure}

\begin{rem}{\rm
We explain our results from the perspective of number density  with Figure \ref{fig:elasticasegmentation}. Consider 10 evenly spaced dots placed within the domain interval \( J \) of a given curve. Each marker corresponds to an equally spaced point in \( J \), but along the curve depicted in the diagram, their spacing follows a pattern determined by the curvature of the curve.

For instance, in range \( A \) of the curve, there are 2 dots; in range \( B \), there are 3 dots; and in range \( C \), there are 5 dots. This demonstrates that the density of dots varies depending on the curvature. The formula in the theorem can be regarded as the continuous counterpart of this discrete phenomenon.
}
\end{rem}

%
%
%
%
%
%

\section{Surface of revolution}
Surface of revolution is an important notion in classical differential geometry. It often appears in several areas of applications of differential geometry such as manufacturing, CG, CV, etc. In this section we deal with surface of revolution with a new coordinate based on the result of plane curve in the previous section, which gives an efficient visualization.

We quickly review basic concepts of surface theory. 
Let $M\subset\R^3$ be a smooth surface parametrized by $f\colon U\to\R^3$ with an open subset $U\subset\R^2$, the curvature coordinate $(u_1,u_2)$. 
Set $\nu$ as the unit normal vector of $M$ and  $\kappa_1, \kappa_2$ as the principal curvatures, where 
$$
\nu_{u_1}=-\kappa_1f_{u_1},\; \nu_{u_2}=-\kappa_2f_{u_2}
$$
hold. $p\in U$ (or $f(p)$) is called {\it a parabolic point} if the Gaussian curvature $K:=\kappa_1\kappa_2$ takes $0$ at $p$. The closure $\mathcal{P}$ of the set of the parabolic points (or its image $f(\mathcal{P})$) is called {\it the parabolic curve}, which often divides the surface into the elliptic region with $K>0$ and hyperbolic region with $K<0$. We call $p\in U$ (or $f(p)$) as a ridge point if $\frac{\partial \kappa_i}{\partial u_i}(p)=0$. The closure $\mathcal{R}$ (or its image $f(\mathcal{R})$) of the set of ridge points are called {\it the ridge curve}.

Let $\gamma\colon I\to \R^2\subset\R^3,\;\gamma(t)=(\gamma_1(t),0,\gamma_2(t))$ be a smooth ($C^\infty$) non-singular curve with an open interval $I$. A surface of revolution is a surface expressed as
$$
f(t,u)=(\gamma_1(t)\cos u,\gamma_1(t)\sin u, \gamma_2(t)).
$$
With this expression, $(t,u)$ is a curvature coordinate of $f$, and the curvature $\kappa$ of $\gamma$ as a plane curve coincides with a principal curvature $\kappa_1$ of $f$ along the $t$-curve (called {\it a profile}); while another principal curvature $\kappa_2$ along the $u$-curve (a cycle) is equal to $\gamma_1(t)$. The surface of revolution $f$ is often visualized with the $t$-curves and $u$-curves. It should be noted that the parametrization of such curvature lines are not unique. Thus we propose one way of efficient visualization of surface of revolution using a fixed curvature coordinate with respect to the parameter $t$.  In particular, we propose to use the tangential angle parameter $\theta$ as the parameter of $\gamma$.

Since $\theta$ is a parameter of a tangential angle of $\gamma$, the above coordinate $(u,\theta)$ gives efficient visualization of the surface, that is, one can easily see how the surface bends along $\theta$ with the interval of  $u$-lines (see Theorem \ref{curveparametrizationthm}). In addition, this parametrization visualize the parabolic curves and ridge curves.

\begin{fact}[\S 2.3 in \cite{Bruce-Giblin-Tari}]
\begin{enumerate}
\item
Let $\gamma(t_i)$ be an inflection of $\gamma$. Then $f(t_i,u)$ is a parabolic point of $f$ for any $u$. 
\item Let $\gamma(t_v)$ be a vertex of $\gamma$. Then $f(t_v,u)$ is a ridge point of $f$ for any $u$. In particular, putting $\kappa_1$ as the principal curvature with $\nu_{t}=\kappa_1f_t$, 
$$
\frac{\partial \kappa_1}{\partial t}(t_v,u)=\frac{d\kappa}{dt}(t_0)=0
$$ 
holds.
\end{enumerate}
\end{fact}

It should be noted that we need divide a surface of revolution by the parabolic curve as the boundary in order to consider the parameter $(\theta,u)$ as a new coordinate. After the suitable stratification based on the parabolic curve, we can take the new coordinate $(\theta,u)$, and have the following Corollary \ref{cor:surfaceofrevo} which is a natural extension of Theorem \ref{curveparametrizationthm}.

\begin{cor}\label{cor:surfaceofrevo}
$$
\frac{\partial }{\partial\theta}\left|\frac{\partial f}{\partial\theta}\right|^2=
-\frac{2}{\kappa_1^4}\frac{\partial\kappa_1}{\partial s},
$$
where $s$ is the arc-length parameter of the profile.
\end{cor}

In the following, we show several figures of surfaces with the new coordinate $(\theta,u)$.

\begin{figure}
\centering{
\includegraphics[width=150mm]{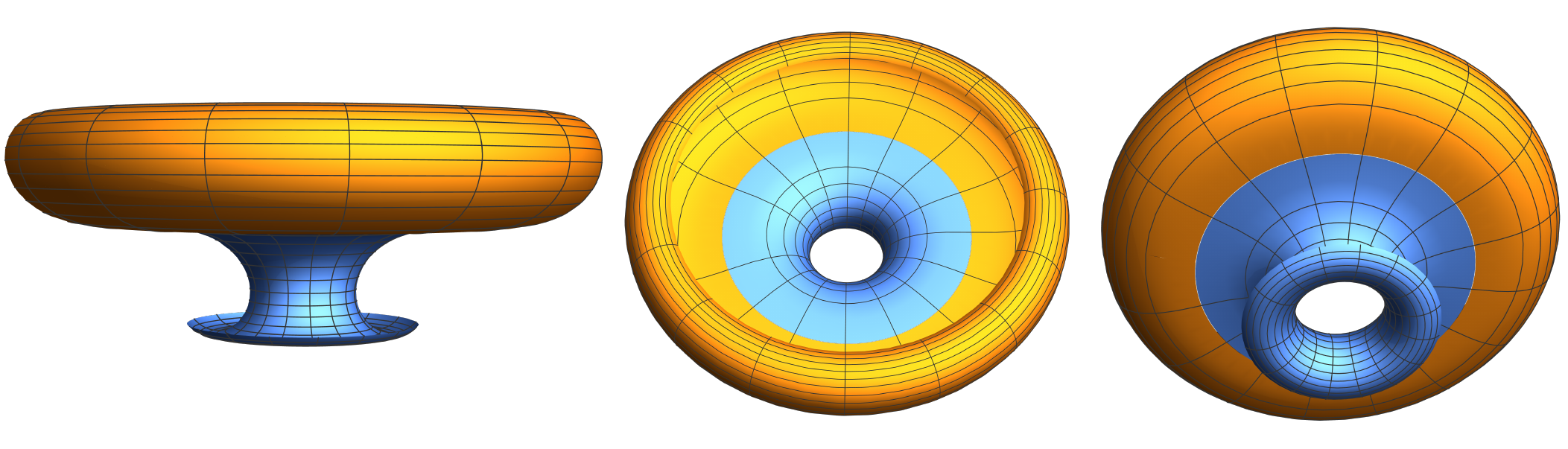}
}
\caption{A surface of revolution with the profile being an Euler spiral. The horizontal circles are plotted along the parameter $\theta$ of the profile. The parabolic curve is the boundary between two regions: the elliptic region colored with orange and the hyperbolic region colored with blue.}
\label{fig:evolution-euler}
\end{figure}

\begin{exam}
Set a plane curve $\gamma\colon \R \to \R^2,\; s\mapsto(\int_0^s\cos\frac{t^2}{2}dt+2,\int_0^s\sin\frac{t^2}{2}dt)$ with the arc-length parameter $s$, which is an Euler spiral. Consider the surface of revolution $f(\theta,u)=(\gamma_1(\theta)\cos u,\gamma_1(\theta)\sin u, \gamma_2(t))$ with $\theta$ being the tangential angle parameter of $\gamma$. Figure \ref{fig:evolution-euler} shows $f$ accompanying $\theta$ and $u$-curves. The blue region indicates the area with negative Gaussian curvature, while the orange region represents positive curvature. The boundary between these regions forms a parabolic curve. This stratification of the surface corresponds to the stratification of the plane curve $\gamma$ based on the sign of its curvature. In particular, the parabolic curve corresponds to inflection points. See Figure \ref{fig:kappa=s} for comparison. 
\end{exam}

\begin{exam}
Consider a plane curve 
$$
\gamma(\theta)=-\frac{\sqrt[3]{2} \left(\sqrt{9 \theta ^2+4}-3 \theta \right)^{2/3}-2 }{2^{2/3} \sqrt[3]{\sqrt{9 \theta ^2+4}-3 \theta }}
(\cos \theta+2, \sin \theta),
$$
whose curvature is $\kappa(s)=1+s^2$ with an arc-length parameter $s$ (see Example \ref{ex:curvewithvertex}). We consider a surface of revolution $f(\theta,u)=(\gamma_1(\theta)\cos u,\gamma_1(\theta)\sin u, \gamma_2(\theta))$.  Figure \ref{fig:evolution-s2+1} shows $f$ accompanying $\theta$ and $u$-curves. The green curve represents the ridge line, which corresponds to the vertex in the plane curve $\gamma$ (see Figure \ref{fig:kappa=s} for comparison). Note that the spacing between the family of $u$-curves is the narrowest on the ridge line.
\begin{figure}
\centering{
\includegraphics[width=150mm]{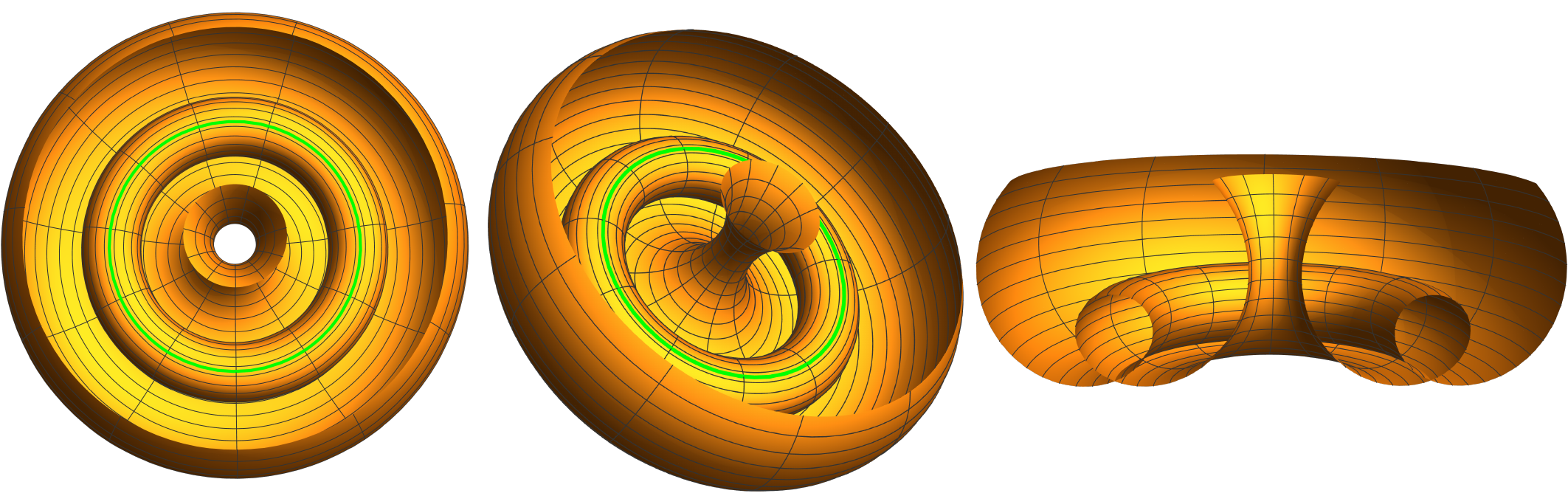}
}
\caption
{A surface of revolution with the profile having a curvature $\kappa(s)=s^2+1$. The horizontal circles are plotted along the parameter $\theta$ of the profile. The green curve means a ridge curve.}
\label{fig:evolution-s2+1}
\end{figure}
\end{exam}

\begin{figure}
\centering{
\includegraphics[width=80mm]{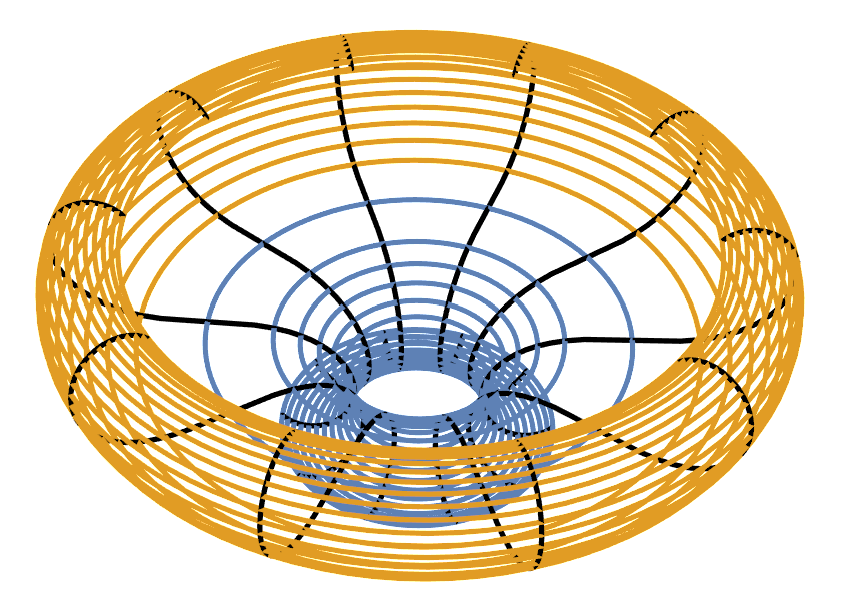}
\includegraphics[width=80mm]{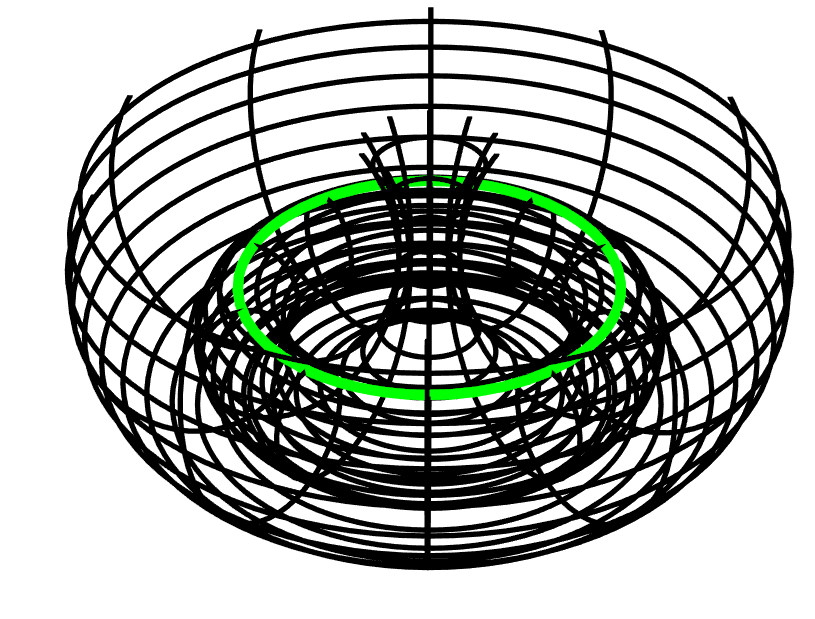}
}
\caption{Surfaces of revolutions with the profiles being an Euler spiral (left) and an elastica (right) drawn by only curvature lines.}
\label{fig:revolutiononlymesh}
\end{figure}

\begin{rem}
{\rm
Figure \ref{fig:revolutiononlymesh} illustrates a surface rendered solely with curvature lines (mesh), where the profile parameterization is based on the tangential angle $\theta$. This representation allows the curvature of the surface to be understood directly from the distribution of the mesh, without relying on shading, lighting effects, or reflections. Such a visualization provides a clear geometric perspective on the bending behavior of the surface.
}
\end{rem}

%
%
%
%
%
\section{Conclusion}\label{sec:conclusion}

Visualizing curvature can be hard when a shape contains local parts with
complicated details in large regions.  In this paper we showed that the
{tangential-angle parameter}~$\theta$ offers a concise remedy.  On planar
curves, sampling at equal increments of~$\theta$ automatically places markers
densely where the curve bends sharply and sparsely where it is almost straight,
so inflection points and vertices emerge without losing sight of their meaning
in the overall context.  On surfaces of revolution, drawing curvature lines at
equal steps of the profile's tangential angle yields a mesh whose spacing
faithfully encodes principal curvature variation and brings parabolic and ridge
curves to the fore.  In every case the resulting graphics provide clear and
consistent curvature visualisations without arbitrary parameter tuning.

The study contributes three main insights.  First, the tangential angle behaves
as an intrinsic parameter that adapts sampling density directly to the local
magnitude of curvature and therefore eliminates manual spacing choices.
Second, a short analytic argument proves that critical points of the function
$\theta \mapsto |d\gamma / d\theta |$ coincide precisely with vertices; this
establishes a transparent link between curvature extrema and marker
distribution.  Third, the very same idea extends naturally to surfaces of
revolution, producing a coherent $(u,\theta)$ coordinate system in which ridge
and parabolic curves appear as intrinsic strata.  A range of examples—
including the elastica, the Euler spiral, a curve with $\kappa(s)=1+s^{2}$, and
two surface models rendered solely with curvature lines—illustrate the
practical impact of these observations.

Several questions remain open.  Although our results are derived for planar
curves and for surfaces generated by revolution, the principle appears to
generalise: if one curvature line parameter on a generic surface is defined as
the integral of a principal curvature,
$u \mapsto \int \kappa_i\,\mathrm{d}s$,
then the resulting coordinate mesh should embed local differential geometric
information into the visualisation.  Demonstrating this extension rigorously
constitutes the first direction for future work.  A second line of inquiry
concerns parallel curves and surfaces, as well as evolutes, whose analysis and
control remain important in both differential geometry and engineering design.
The parameter introduced here promises new insight into how curvature evolves
under offset operations and may ultimately help to guide geometric tolerances
in practical applications.


\begin{thebibliography}{nn}
{\small

%




\bibitem{Bruce-Giblin} J. W. Bruce and P. J. Giblin, \newblock{Curves and Singularities. A geometrical introduction to singularity theory. Second edition}. \newblock{Cambridge University Press}, Cambridge, 1992.

\bibitem{Porteous} I. R. Porteous, \newblock{Geometric Differentiation. For the Intelligence of Curves and Surfaces}. \newblock{Cambridge University Press}, Cambridge, 2001.


\bibitem{Bruce-Giblin-Tari} J. W. Bruce, P. J. Giblin and F. Tari, \newblock{Ridges, crests and sub-parabolic lines of evolving surfaces.} \newblock{Int. J. Comput. Vision} {18} (1996), 195--210.


\bibitem{Cipolla-Giblin}
R. Cipolla and P. Giblin,
\newblock{Visual motion of curves and surfaces},
\newblock{Cambridge Univ. Press.} 2000.

\bibitem{Damon-Giblin-Haslinger}
J. Damon, P. Giblin and G. Haslinger,
\newblock{Local features in natural images via singularity theory},
\newblock{Lecture Notes in Math.  2165, Springer.} 2016.

\bibitem{Hertzmann2000}
A.~Hertzmann and D.~Zorin.
\newblock Illustrating smooth surfaces.
\newblock In {Proceedings of SIGGRAPH 2000}, pages 517--526, 2000.

\bibitem{Dong2006}
S.~Dong, P.~T. Bremer, M.~Garland, V.~Pascucci, and J.~C. Hart.
\newblock Spectral surface quadrangulation.
\newblock {ACM Transactions on Graphics}, 25(3):1057--1066, 2006.

\bibitem{Palacios2007}
J.~Palacios and E.~Zhang.
\newblock Rotational symmetry field design on surfaces.
\newblock {ACM Transactions on Graphics}, 26(3):55:1--55:8, 2007.

\bibitem{Iarussi2015}
E. Iarussi, D. Bommes, and A. Bousseau, 
\newblock{BendFields: Regularized Curvature Fields from Rough Concept Sketches.}
\newblock{ACM Transactions on Graphics}, 
vol. 34, no. 3, article 24, 2015, 16 pages.
\newblock{https://doi.org/10.1145/2710026}


\bibitem{IRFT} S. Izumiya, M. C. Romero Fuster, M. A. S. Ruas and F. Tari,
\newblock{Differential Geometry from a Singularity Theory Viewpoint}.
\newblock{World Scientific Pub. Co Inc.} 2015.

\bibitem{Takezawa2019}
M. Takezawa, K. Matsuo, and T. Maekawa, 
\newblock{Control of lines of curvature for plate forming in shipbuilding.}
\newblock{Computer Aided Geometric Design}, 
vol. 75, 2019, article 101785. 
\newblock{https://doi.org/10.1016/j.cagd.2019.101785}


%


%
%
\bibitem{Koenderink}
J. J. Koenderink,
\newblock{Solid shape},
\newblock{MIT Press Series in Artificial Intelligence. MIT Press.} 1990.
%
\bibitem{Koenderink1992}
J.~J. Koenderink and A.~J. van Doorn.
\newblock Surface shape and curvature scales.
\newblock {Image and Vision Computing}, 10(8):557--564, 1992.



\bibitem{KD1998}
J. J. Koenderink and A. Doon, Shape from Chebyshev nets, Proc. ECCV’98 (1998), 215--225.



}
\end{thebibliography}
\end{document}